\documentclass[a4paper,12pt]{article}

\newtheorem{definition}{Definition}
\newtheorem{assumption}{Assumption}
\newtheorem{remark}{Remark}
\newtheorem{problem}{Problem}

\usepackage{graphicx}
\usepackage{amsmath,amsfonts,amssymb}
\usepackage{theorem}
\usepackage{xspace}
\usepackage{xifthen}
\usepackage{graphicx}
\newcounter{probCounter}
\renewenvironment{problem}[1]
{\refstepcounter{probCounter}~\\[1.5ex]\noindent\textbf{Problem~\arabic{probCounter}\ (#1):\ }\begingroup\it}
{\endgroup~\\[1.5ex]}
\newcounter{defCounter}
\renewenvironment{definition}[1]
{\refstepcounter{defCounter}~\\[1.5ex]\noindent\textbf{Definition~\arabic{defCounter}\ (#1):\ }\begingroup\it}
{\endgroup~\\[1.5ex]}
\newcounter{assuCounter}
\renewenvironment{assumption}[1]
{\refstepcounter{assuCounter}~\\[1.5ex]\noindent\textbf{Assumption~\arabic{assuCounter}\ (#1):\ }\begingroup\it}
{\endgroup~\\[1.5ex]}
\newcounter{rmkCounter}

\newcommand{\iset}[2]{\ensuremath{#1,\ldots,#2}}

\newcommand{\mom}[2][]{%
\ifthenelse{\equal{#1}{}}%
{\ensuremath{\nu\left(#2\right)}}%
{\ensuremath{\nu^{(#1)}\left(#2\right)}}}%
\newcommand{\momv}[2][]{%
\ifthenelse{\equal{#1}{}}%
{\ensuremath{{\nu}\left(#2\right)}}%
{\ensuremath{{\nu}^{(#1)}\left(#2\right)}}}%
\newcommand{\ie}{i.e.\xspace}
\newcommand{\cf}{cf.\xspace}
\newcommand{\eg}{e.g.\xspace}

\newcommand{\transp}{^\mathsf{T}}
\newcommand{\coloneqq}{:=}

\title{\bf 
  Probabilistic and Set-based Model Invalidation and Estimation using LMIs}

\begin{document}

\author{Stefan Streif$^1$, Didier Henrion$^{2,3,4}$ and Rolf Findeisen$^1$}

\footnotetext[1]{Institute for Automation Engineering,
Otto-von-Guericke-University~Magdeburg, Germany.}
\footnotetext[2]{CNRS, LAAS, University of Toulouse, France.}
\footnotetext[3]{Faculty of Electrical Engineering, Czech Technical University in Prague,
Czech Republic}

\date{Draft of \today}

\maketitle

\begin{abstract}
Probabilistic and set-based methods are two approaches for model \mbox{(in)}validation, parameter and state estimation.
Both classes of methods use different types of data, \ie deterministic or probabilistic data, which allow different statements and applications.
Ideally, however, all available data should be used in estimation and model invalidation methods. 
This paper presents an estimation and model \mbox{(in)}validation framework combining set-based and probabilistically uncertain data for polynomial continuous-time systems.
In particular, uncertain data on the moments and the support is used without the need to make explicit assumptions on the type of probability densities.
The paper derives pointwise-in-time outer approximations of the moments of the probability densities associated with the states and parameters of the system.
These approximations can be interpreted as guaranteed confidence intervals for the moment estimates.
Furthermore, guaranteed bounds on the probability masses on subsets are derived and allow an estimation of the unknown probability densities. 
To calculate the estimates, the dynamics of the probability densities of the state trajectories are found by occupation measures of the nonlinear dynamics.
This allows the construction of an infinite-dimensional linear program which incorporates the set- and moment-based data. This linear program is relaxed by a hierarchy of LMI problems providing, as shown elsewhere, an almost uniformly convergent sequence of outer approximations of the estimated sets.
The approach is demonstrated with numerical examples.
\end{abstract}

\section{Introduction}

State and parameter estimation as well as model \mbox{(in)}validation are frequently used in many applications with different requirements.
Often different kind of data, \ie probabilistic or deterministic, are available providing different insights into the process.
The determination of the set of all parameters consistent with measurements or proving nonexistence of consistent parameters allows for guaranteed model invalidation \eg in fault detection and isolation \cite{Savchenko_etAl_2011_SB_FDI_PolynSys_Hybrid}.
For processes such as crystallization \cite{Fujiwara_etAl_2005_Control_crystallization,Marchisio_2007_Quadrature_MoM} or cell populations \cite{Zhu_etAl_2000_MPC_CellPop_PopBalance}, however, not only the estimation of consistent parameters, but also the estimation of probability densities is important. Uncertain estimates of the moments of a probability densities can often be determined based on samples \eg from a number of similar experiments (see references in \cite{Delage_Ye_2010_ProbROptim_MomUnc}).
Ideally, an estimation method should allow all available information and associated uncertainties to be taken into account, \ie both deterministic and probabilistic uncertainties.

Set-based estimation approaches employing a bounded-error uncertainty description provide guaranteed yes/no answers on model invalidity and outer or inner approximations of consistent parameter sets.
Typically, probability densities and data are not used in set-based estimation approaches and only the support of the uncertain variables is considered.
Probabilistic estimation methods take data on the moments of probability densities into account and allow the estimation of probability densities and therefore probabilistic statements on model validity.

This work presents initial steps toward a combination of set- and moment-based data into a consistent framework for estimation and model \mbox{(in)}validation of nonlinear (polynomial), continuous-time systems.
To the best of our knowledge, the presented framework has not been considered in this context before. However, numerous other methods addressing similar problems exist in the literature and several of them will be reviewed briefly below.
Due to space limitations, an exhaustive literature review is beyond the scope of this paper.
Further references and a literature review on model invalidation methods and set-based estimation methods can be found, \eg, in \cite{Streif_etAl_2013_CDC__ContParamEstim,Halder_Bhattacharya_2011_ProbMInv}.

Several methods for set-based estimation and model invalidation exist, such as set-member\-ship methods \cite{Milanese_Vicino_1991_OptimalEstim_DynSys_SetMemb}, relaxation based approaches \cite{Streif_etAl_2013_CDC__ContParamEstim,Cerone_etAl_2012__Set-membership_identification_convex_relaxation} (and references therein), robust and $\mathcal{H}_\infty$ control where model validation methods have been studied extensively assuming structured norm-bounded uncertainty in a usually linear systems setting (see \cite{Zho_Doyle_Glover_1995_RContr_OC,Smith_Doyle_1992_MVal_RContr_Ident} and references therein).
Funnels for consistent state trajectories (but not moments) were proposed in \cite{Tobenkin_etAl_2011_Funnels_Traj_SOS}.
Barrier certificates approaches using sum-of-squares restrictions \cite{Prajna_2006_Automatica__Barrier_certificate_model_validation} are similar to the LMI approach proposed here; to the best of our knowledge, however, the existence of barrier certificates is not guaranteed.
The work by \cite{Dabbene_etAl_2012_PartI_Prob_PE,Dabbene_etAl_2012_PartII_Prob_PE} does not explicitly consider dynamical systems, but it provide methods and randomized and deterministic algorithms to use statistical assumptions within set-based estimation. By this, the set-based estimates can be improved at the expense of an allowed probabilistic risk. 

\cite{Singh_Hespanha_2011_ApproxMomDyn} consider the stochastic dynamics of the state variables and the state density is modeled via its moments.
In contrast, the initial setup considered in this paper is not stochastic in the sense that the time evolution of the moments is modeled via the Liouville equation.
This equation reformulates nonlinear ordinary differential equation (ODE), with the help of occupation measures, into a linear partial differential equation to describe the time evolution of the density of the solution of the ODE.
This does not lead to a moment enclosure problem as in \cite{Singh_Hespanha_2011_ApproxMomDyn} and allows -- due to the convex relaxations employed in this work -- deriving outer approximations of the moments.

Uncertainties in the data, including data on the moments, have been considered in robust optimization, see \cite{Delage_Ye_2010_ProbROptim_MomUnc}. A common approach is to transform risk measures or probabilistic uncertainties into deterministic uncertainty sets \cite{Natarajan_etAl_2009_Constructing_risk_measures_from_uncertainty_sets}. 
Bounds on moments using convex optimization and LMIs have also been presented by \cite{Vandenberghe_Boyd_1999_Applications_SDP,Bertsimas_Caramanis_2006_Bounds_LinPDE_SDP,Bertsimas_etAl_2000_MoM_SDP}.

An approach to probabilistic model (in)validation using probability metrics has been presented by \cite{Halder_Bhattacharya_2012_ProbMInv_Wasserstein_metric,Halder_Bhattacharya_2011_ProbMInv}, which can be used to compare probabilistic model predictions with data. As in our work, Liouville's equation is used to propagate the uncertainties. However, convex relaxations allowing the outer approximations of moments were not considered.

Many probabilistic estimation methods require (implicitly or explicitly) assumptions on the type of probability densities of the estimates.
Parameter and state estimation from a system identification perspective \cite{Ljung_1999_SysId} aims at validations through a statistical correlation analysis of the residuals.
Randomized algorithms employ randomness in the simulation of particles and use a Bayes approach for estimation. Well-known techniques are Markov-Chain-Monte-Carlo approaches and particle filters for optimal and suboptimal Bayesian algorithms for nonlinear/non-Gaussian tracking problems. For a review see  \cite{Arulampalam_etal_2002_Tutorial_PF_NonlSys_Tracking}. 

Polynomial Chaos has been used for parameter estimation \cite{Ghanem_etAl_ProbMInv_PCE,Kim_etal_2013_CSM__PCE_review}, but usually requires knowledge of the probability densities of the uncertainties. These methods allow the transformation of the stochasticity in the parameters or initial conditions by a finite series of orthogonal polynomials into deterministic equations enabling an easy analysis of the moments of the state trajectory obtained by numeric integration.
%


The initial setup of the presented framework is set-based and has been used in \cite{Streif_etAl_2013_CDC__ContParamEstim} for the estimation of consistent parameter sets and model invalidation.
In this work, the framework is extended such that probabilistic, moment and set-based uncertainties and data can be considered. 
This allows the outer approximation of the moments and of the support of the probability densities of the parameters, initial conditions or state trajectories pointwise-in-time, without explicit assumptions on the underlying probability density and without numerical integration or approximation by discrete-time models.
Furthermore, bounds on the probability mass on a subset can be given, which can be used to approximate the shape of the probability densities of the parameters or for probabilistic model validation.
To reach these goals, converging hierarchies of LMI problems are derived by reformulating the nonlinear dynamics in terms of occupation measures as described in \cite{Henrion_Korda_2013_ROA_occupation_measure,Streif_etAl_2013_CDC__ContParamEstim,Lasserre_etAl_2008_SIAM__Nonlin_OC_OccupMeas_LMIs}.

The remainder of this contribution is structured as follows.
Sec.~\ref{sec:setup} defines the dynamical model and the constraints and uncertainties imposed by the probabilistic and set-based measurements and additional a-priori information, and it states the considered problems.
Sec.~\ref{sec:mom_estim} formally introduces the framework of occupation measures and derives outer approximations of pointwise-in-time moments. Those outer approximations define the imprecision or confidence in the moment estimates.
Furthermore, lower and upper bounds on the probability mass on a subset are derived.
Small example are used in Sec.~\ref{sec:mom_estim} to illustrate the results.
The approach and future research directions are discussed in Sec.~\ref{sec:conclusion}.

\subsubsection*{Notation.}
$n_x$ denotes dimension of a vector $x$.
Subscripts of vectors and matrices denote the corresponding row and column elements.
Sets and function spaces are denoted by capital, calligraphic letters, \eg $\mathcal{X}$.
$\mathcal{M}(\mathcal{A})$ denotes the set of finite Borel measures supported on the set $\mathcal A$, which can be interpreted as
elements of the dual space $\mathcal C(\mathcal{A})'$, \ie as bounded linear functionals acting on the set of continuous
functions $\mathcal C(\mathcal{A})$. 
${P}(\mathcal{A})$ denotes the set of probability measures on $\mathcal A$, \ie those measures $\mu$ of $\mathcal{M}(\mathcal{A})$ which are nonnegative and normalized to $\mu(\mathcal{A})=1$.
The $m^\text{th}$ moment ($m\in\mathbb{N}$) of $\mu(x) \in P(\mathcal{X})$ supported on the set $\mathcal{X}$ is denoted by
$\nu^{(m)}(x)
\coloneqq \int_{\mathcal{X}} x^m \mu(dx)$.
All moments up to degree $d$ are denoted by $\nu^{(\le d)}$.
We define the probability mass on a set $\widehat{\mathcal{X}}\subseteq\mathcal{X}$ by $F(\widehat{\mathcal{X}}) \coloneqq \int_{\widehat{\mathcal{X}}} \mu(dx) \le 1$.
Pointwise-in-time probability measures are denoted by $\mu_k(dx)$ (and their moments by $\nu^{(m)}_k$), and occupation measures on the time-interval $[t_k,t_{k+1}]$ are denoted by $\mu_{k,k+1}(dt,dx)$ (and their moments by $\nu^{(m)}_{k,k+1}$).

\section{Problem Setup}
\label{sec:setup}

\subsection{Dynamical System}

Consider continuous-time, nonlinear systems of the following form
\begin{align}
\label{eq:system}
    \dot{x}(t) &= f\bigl(t,x(t)\bigr), \quad x(0) = x_0.
\end{align}
The states and initial conditions are denoted by $x \in \mathbb{R}^{n_{x}}$ and $x_0 \in \mathbb{R}^{n_x}$, respectively.
Time is denoted by $t$ and is restricted to the interval $[0,1]$, which can be achieved by a suitable time-scaling of the dynamics.

We assume the vector fields $f$ to be polynomial maps.
Time-invariant variables, \ie parameters, can be accounted for in \eqref{eq:system} by adding trivial dynamics $\dot{x}_i = 0$ for some $i \in \{1,\ldots, n_x\}$.
Note that we do not explicitly consider output maps. 


We assume set-based and moment-based uncertainty descriptions of the variables, \ie parameters, initial conditions and measurements.
With these data and using the nonlinear continuous-time model \eqref{eq:system}, outer approximations of the consistent moments and the support of the variables are estimated pointwise-in-time. 
In Sec.~\ref{sec:problems} we formally define consistency and state the considered problems.

\subsection{Uncertainty Descriptions}
\label{sec:uncertainties}

The set-based uncertainties define the support of the variables, whereas the moment-based probabilistic uncertainties constrain the probability measures. 
We assume that set-based uncertainties on the variables $x(t)$ at time $t_k$ are given by $m_x$ polynomial inequalities:
\[
    \label{eq:data_z}
    \mathcal{X}_k \coloneqq \bigl\{ x(t_k) :\ g_{i}(t_k,x(t_k)) \geq 0, i = \iset{1}{m_x} \bigr\} \subset
    \mathbb{R}^{n_x}.
\]
Eq.~\eqref{eq:data_z} can be used to represent constraints on the initial conditions ${x_0}$ (including the parameters) and measurements of states (or polynomial combinations thereof) at time points $t_k$, $k=\iset{1}{m_t}$.

In a similar manner to \eqref{eq:data_z}, we define a set $\mathcal{X}$ such that $x(t) \in \mathcal{X}\subset\mathbb{R}^{n_x}, \forall t \in [0,1]$ and such that $\mathcal{X}_k \subseteq \mathcal{X}$.

Constraints on the moments (up to order $d$) of the probability measure $\mu_k$ at time $t_k$ are given by $m_\nu$ polynomial inequalities:
\[
    \label{eq:data_z_mom}
    \mathcal{M}_k \coloneqq \Bigl\{ \nu^{(\le d)}_k :\ h_{i}\bigl(\nu^{(\le d)}_k\bigr) \geq 0, i = \iset{1}{m_\nu} \Bigr\} \subset \mathbb{R}^{n_\nu}. 
\]
As above for the set-based uncertainties, Eq.~\eqref{eq:data_z_mom} can be used to represent constraints on the moments of the initial conditions ${x_0}$ (including the parameters), the states or the measurements pointwise-in-time.

\subsection{Problem Statement}
\label{sec:problems}

\begin{figure}[t]
\centering
\includegraphics[width=0.9\linewidth]{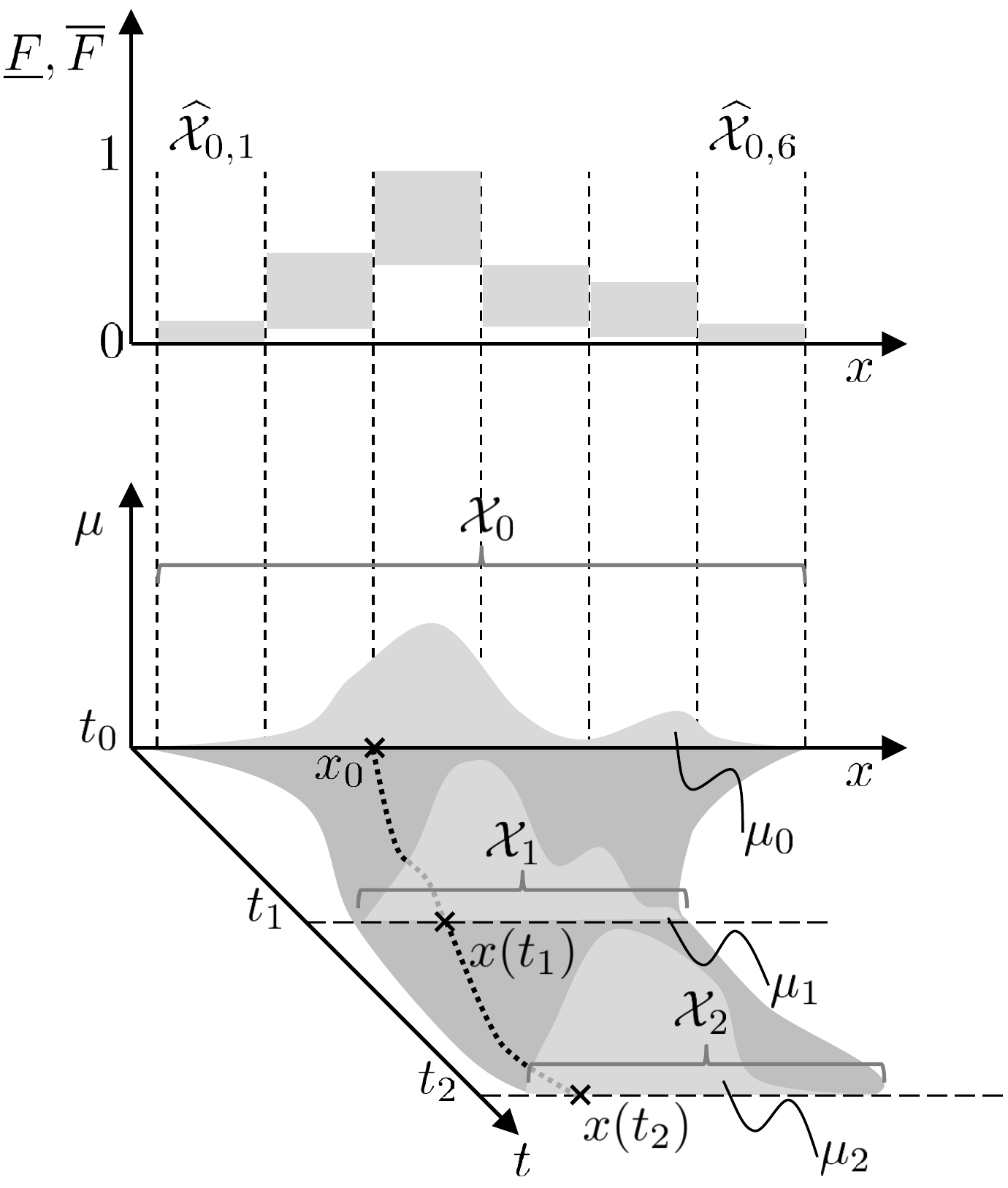}
\caption{Illustration of notation, the problems and the approach. In the lower part of the figure, the pointwise-in-time probability densities $\mu_k$ (and their support $\mathcal{X}_k$) are shown by the light-gray shaded areas at $t_k$, $k = 0,1,2$. The dotted line shows a sample trajectory starting at $x_0$ and ending at $x(t_2)$, and the dark gray shaded area shows the set of possible trajectories. In the upper part, the lower bounds $\underline{F}$ and upper bounds $\overline{F}$ on the probability mass over the subsets $\widehat{\mathcal{X}}_{0,1}, \ldots, \widehat{\mathcal{X}}_{0,6} \subset \mathcal{X}_0$ are shown by the gray boxes.}
\label{fig:problems}
\end{figure}

For general purposes, we define the following consistency set using the nonlinear dynamics \eqref{eq:system}, the set-based \eqref{eq:data_z} and the probabilistic uncertainties \eqref{eq:data_z_mom}:
\begin{definition}{Consistency set}
\label{def:consistent:Z_zeta}
The set of consistent values of $x(t_k)$ and $\nu^{(\le d)}_k$ at time-point $t_k$ is given by:
\begin{align}
\label{eq:consistent:Z_zeta}
\mathcal{C}^\ast_k \coloneqq\ &\Bigl\{\ \bigl(x(t_k),\nu^{(\le d)}_k\bigr) : \exists\ \bigl(x(t), \mu(t, x(t))\bigr)\ \forall t \in [0,1] \nonumber\\
&\quad     
  \text{such that\ }\quad x(t_l) \in \mathcal{X}_l, l = 0, \ldots, m_t, \text{\ \ and} \nonumber\\
&\quad     
  \phantom{\text{such that\ }}\quad \nu^{(\le d)}_l \in \mathcal{M}_l, l = 0, \ldots, m_t, \text{\ and}\nonumber\\
&\quad     
  \phantom{\text{such that\ }} x(t)\!=\!x_0\!+\!\int_0^{t}\!\!f\bigl(t,x(t)\bigr) dt \in \mathcal{X}\ \forall t\!\in\![0,1]\Bigr\} \subset \mathbb{R}^{n_x+n_\nu}.
\end{align}
In particular, the set of consistent initial conditions is given by the orthogonal projection of $\mathcal{C}^\ast_0$ onto $x_0$, \ie $\perp_{x_0} \mathcal{C}^\ast_0$.
Furthermore, the consistent pointwise-in-time moments of $x(t_k)$ are given by the projection $\perp_{\nu^{(\le d)}_k} \mathcal{C}^\ast_k$.
\end{definition}
Based on Def.~\ref{def:consistent:Z_zeta}, we state the problems considered in this paper. See also Fig.~\ref{fig:problems} for illustrations.

The first problem aims to estimate the consistent sets of moments of the parameters and states based on the measurement data and uncertainties \eqref{eq:data_z}--\eqref{eq:data_z_mom} under consideration of the dynamics \eqref{eq:system}.
\begin{problem}{Pointwise-in-time estimation of moments}
    \label{prob:mom_estim}%
    Find consistent pointwise-in-time moments $\mathcal{C}^\ast_0$ of the consistent initial conditions (including parameters) and $\mathcal{C}^\ast_k, t_k \in [0,1]$ of the states.
\end{problem}
Due to the possible uncertainties in the data, the results are sets of consistent moments, where the size of the set reflects the confidence of the moment estimates.
Note that the pointwise-in-time estimation of the consistent sets of the initial conditions and parameters has been treated in \cite{Streif_etAl_2013_CDC__ContParamEstim}; however, without probabilistic data or moment constraints.

Set-based approaches can be used to provide model invalidity certificates, \ie yes/no answers whether a model is invalid or not. In case a model is valid, \ie there exist consistent parameter values, then it is often of interest to quantify probabilistic model validity which is related to the shape of probability density \cite{Halder_Bhattacharya_2011_ProbMInv,Halder_Bhattacharya_2012_ProbMInv_Wasserstein_metric}. We therefore consider the following problem:
\begin{problem}{Estimation of probability measure}
    \label{prob:MInv}
    Determine lower and upper bounds ($\underline{F}$ and $\overline{F}$, respectively) on the probability mass $F$ over the set $\widehat{\mathcal{X}}_{0} \subseteq \mathcal{X}_0$.
\end{problem}
Prob.~\ref{prob:MInv} allows probabilistic model validation under consideration of set-based and moment-based measurement data: if $\overline{F}$ is small, then the (cumulative) probability over the subset $\widehat{\mathcal{X}}_{0}$ for the constraints \eqref{eq:data_z} and \eqref{eq:data_z_mom} to be satisfied is small. The comparison of the probability mass of different subsets $\widehat{\mathcal{X}}_{0,1}, \ldots, \widehat{\mathcal{X}}_{0,m_s}$ then allows quantifying the likeliness that a random parameter sample from this set satisfies the moment- and set-based constraints.
The subset $\widehat{\mathcal{X}}_{0}$ could be obtained from a partitioning of the set $\mathcal{X}_0$ as shown in Fig.~\ref{fig:problems}. Furthermore, it allows the approximation of the probability density, \cf Fig.~\ref{fig:problems}.

Note that in Prob.~\ref{prob:MInv} bounds on the $0^\text{th}$-order moment (\ie, the probability mass) are derived for a given subset, whereas in Prob.~\ref{prob:mom_estim} the consistent set of state, parameter and moments are determined based on measurements.

In Prob.~\ref{prob:mom_estim} and \ref{prob:MInv}, the desired estimates are difficult to determine due to the nonlinearities of the dynamics and other constraints, which results in nonconvexity of the sets.
In the next section, we use occupation measures and convex relaxations to account for the time-continuous dynamics and this allows the efficient computation of outer approximations of the estimates.

Note that in this work, we make the following assumption:
\begin{assumption}{Bounded support}
The probability measure $\mu(t,x)$ is supported  on $[0,1]\times\mathcal{X}$, \ie $\mu(t,x) \in P([0,1]\times \mathcal{X})$.
\end{assumption}

\section{Outer Approximation of Moments and Estimation of Probability Mass}
\label{sec:mom_estim}

In the first subsection an intuitive and condensed explanation of the method and concepts of the underlying mathematical framework is presented.
All mathematical details and proofs are found in \cite{Henrion_Korda_2013_ROA_occupation_measure,Streif_etAl_2013_CDC__ContParamEstim,Lasserre_etAl_2008_SIAM__Nonlin_OC_OccupMeas_LMIs} and in the references therein.
In the subsequent subsections, Prob.~1 and 2 are addressed and the results are illustrated for simple examples.

\subsection{Constrained Uncertainty Propagation in Nonlinear Continuous-Time Systems}

The approach can be summarized as follows. The original deterministic nonlinear dynamics, given by \eqref{eq:system}, is reformulated in the space of nonnegative measures.
To do so, so-called occupation measures $\mu_{0,1}(dt,dx)$ are introduced.
These probability measures encode the nonlinear dynamics on the time interval $[0,1]$ and allow the consideration of uncertainties by averaging over space and time which also eliminates the nonlinear expressions.
Furthermore, they are linked by pointwise-in-time probability measures $\mu_0(dx)$ and $\mu_1(dx)$ via Liouville's equation:
\begin{multline}
\label{eq:Liouville}
\int_\mathcal{X} \int_0^1 \left( \frac{\partial v}{\partial t}(t,x) + \nabla_x v(t,x) f(t,x)\right) \mu_{0,1}(dt,dx)
\\= \int_\mathcal{X} v(1,x)\mu_{1}(dx) - \int_\mathcal{X} v(0,x)\mu_0(dx)
\end{multline}
which allows deriving linear constraints as shown below.
Here, Liouville's equation is written in its variational form, for all monomials $v(t,x)=t^{\alpha}x^{\beta}$.
Note that the monomials are a dense basis for the set of continuous functions on compact sets.
Liouville's equation is a linear partial differential equation which accounts for the time-evolution of probability measures ruled by the nonlinear continuous-time dynamics.
Note that if the initial conditions are certain (\ie points for fixed $t$ and hence a Dirac probability distribution), Liouville's equation describes the evolution of the test function $v(t,x)$ along the trajectory.

This reformulation in terms of occupation measures results in a system of a finite number of linear constraints in an infinite-dimensional space,
which links all the moments $\nu^{(\infty)}_{0,1}$ of the occupation measures $\mu_{0,1}(dt,dx)$ and the moments $\nu^{(\infty)}_{j}$ of the measure $\mu_j(dx)$ ($j=0,1$), respectively.
In case of the probability measure $\mu_0(dx)$, the moments are related to the test functions $v(t,x)$ by
\begin{align}
\label{eq:Liouville:link_moments_testFuncs}
\nu^{(m)}_0 \coloneqq \int_{\mathcal{X}} x_0^m \mu_0(dx) \quad \text{with}\quad v(t,x) = t^0 x^m,
\end{align}
and similar for the other probability measures.
Let us denote the resulting infinite-dimensional linear system of equations by $\widetilde{A}(\nu^{(\infty)}_{0,1},\nu^{(\infty)}_{0},\nu^{(\infty)}_{1})=\widetilde{b}$.
The obtained set of equations can be extended by equations constraining the support of the probability measures (\cf set-based uncertainties \eqref{eq:data_z}), and by enforcing constraints on the moments of the probability measures thereby including data on the moments (\cf probabilistic uncertainties \eqref{eq:data_z_mom}).
For a measure to be supported on a semi-algebraic set, the infinite-dimensional vectors of moments $\nu^{(\infty)}_0$, $\nu^{(\infty)}_1$ and $\nu^{(\infty)}_{0,1}$ have to fulfill necessary and sufficient conditions.
These are formulated by infinite-dimensional convex linear matrix inequalities (LMI) called moment matrices that we 
express by $M(\nu^{(\infty)}_0)\succeq 0$ to guarantee positivity of the measures.
Furthermore, the matrices necessary to account for the support are called localizing matrices and are denoted by $L(g_i\nu^{(\infty)}_0) \succeq 0$, where the $g_i$ are polynomials.
The matrices $M$ are symmetric, square matrices whose rows and columns are indexed by the monomials $v(t,x)$ (\cf Eq.~\eqref{eq:Liouville:link_moments_testFuncs}).
The matrices $L$ are similar to $M$, but their elements are additionally multiplied by the polynomial $g_i$.

The resulting infinite-dimensional decision problem can be solved by a converging hierarchy of semidefinite relaxations or sum-of-squares restrictions.
This then allows determining an outer approximation of the moments and initial conditions (\cf Prob.~\ref{prob:mom_estim} and 2).
If we truncate the infinite sequence of moments to moments of degree up to $d$,
we obtain a finite-dimensional linear system of equations that we denote by
$A(\nu^{(\le d)}_{0,1},\nu^{(\le d)}_{0},\nu^{(\le d)}_{1})=b$, as well as truncated finite-dimensional LMI constraints $M(\nu^{(\le d)}_0)\succeq 0$ and $L(g_i\nu^{(\le d_{g_i})}_0)\succeq 0$, where $d_{g_i} \coloneqq d-\text{degr}(g_i)$.
Note that these constraints are necessary for the corresponding measures to be supported on a semi-algebraic set.
By construction, minimization (resp. maximization) of an entry of the vectors $\nu^{(\le d)}_{0}$ or $\nu^{(\le d)}_{1}$ on the resulting
finite-dimensional convex set yields a valid lower (resp. upper) bound on the corresponding moment. 
When increasing the truncation degree $d$ (also called relaxation order), we obtain a monotonically
nondecreasing (resp. nonincreasing) sequence of lower (resp. upper) bounds that converge
to the exact value of the moment consistent with the uncertain set-based and moment-based measurement data.

To address Prob.~1 and 2 and to recover local, non-averaged information in time and space, a discretization is needed, which still allows guaranteed outer approximations. This is explained in the following subsections.

\subsection{Estimation of Moments Pointwise-in-time}

As in \cite{Streif_etAl_2013_CDC__ContParamEstim} we split the global occupation measure time-wise into local occupation measures $\mu_{k,k+1}(dx)$ corresponding to time intervals $[t_k,t_{k+1}]$, $k=0,1,\ldots,m_t-1$.
For the upper bound on a moment $\widetilde{\nu}\in \{\nu^{(\le d)}_{0},\ldots,\nu^{(\le d)}_{m_t}\}$, this then gives the following moment relaxation of order $d$:
\begin{equation}
\label{eq:prob1:mom_relax}
\setlength{\arraycolsep}{1pt}
\begin{array}{lll}
\text{maximize} & \widetilde{\nu} \\
\text{subject to} &\ A\Bigl(\nu^{(\le d)}_{k,k+1},\nu^{(\le d)}_{k},\nu^{(\le d)}_{k+1}\Bigr) = b_{k,k+1}, & k = 0, \ldots, m_t-1,\\
& M\Bigl(\nu^{(\le d)}_{k,k+1}\Bigr) \succeq 0, & k = 0, \ldots, m_t-1,\\
& M\Bigl(\nu^{(\le d)}_{k}\Bigr) \succeq 0, & k = 0, \ldots, m_t,\\
& L\Bigl(g_i\nu^{(\le d_{g,i})}_{k}\Bigr) \succeq 0, & i = 0,\ldots m_x, k = 0, \ldots, m_t,\\
& L\Bigl(h_i\nu^{(\le d_{h,i})}_{k}\Bigr) \succeq 0, & i = 0,\ldots m_\nu,  k =  0, \ldots, m_t,
\end{array}
\end{equation}
where $g_i$, $i = 1, \ldots, m_x$ and $h_i$, $i = 1,\ldots, m_\nu$ are polynomials from \eqref{eq:data_z} and \eqref{eq:data_z_mom}, respectively.
The matrices $M$ and $L$ are the truncated moment and localizing matrices, respectively, depending linearly on moment vectors of the respective degrees.

\subsubsection{Example 1.}
\label{sec:mom_estim:example1}

Consider the following bilinear example
\begin{subequations}
\begin{align}
\label{eq:mom_estim:example1}
\dot{x}_1 &= -x_1x_2  \\
\dot{x}_2 &= 0
\end{align}
\end{subequations}
with certain, fixed initial conditions $x_1(0)=0.5$ (represented by a Dirac probability measure), and probabilistically uncertain parameter $x_2$ uniformly distributed on $[0,1]$.
This example was chosen because an analytic solution is available which allows an easy comparison with the computational results. The analytic solution is:
\begin{subequations}
\label{eq:mom_estim:example1:solution}
\begin{align}
x_1(t) &= x_1(t_0)\exp(-x_2t), \\
\mom[1]{x_1(t)} &= x_1(t_0)\frac{1 - \exp(-t)}{t},\\
\mom[2]{x_1(t)} &= x_1(t_0)^2\frac{1 - \exp(-2t)}{2t}.
\end{align}
\end{subequations}

\begin{figure}[h]
\centering
\includegraphics[width=0.8\linewidth]{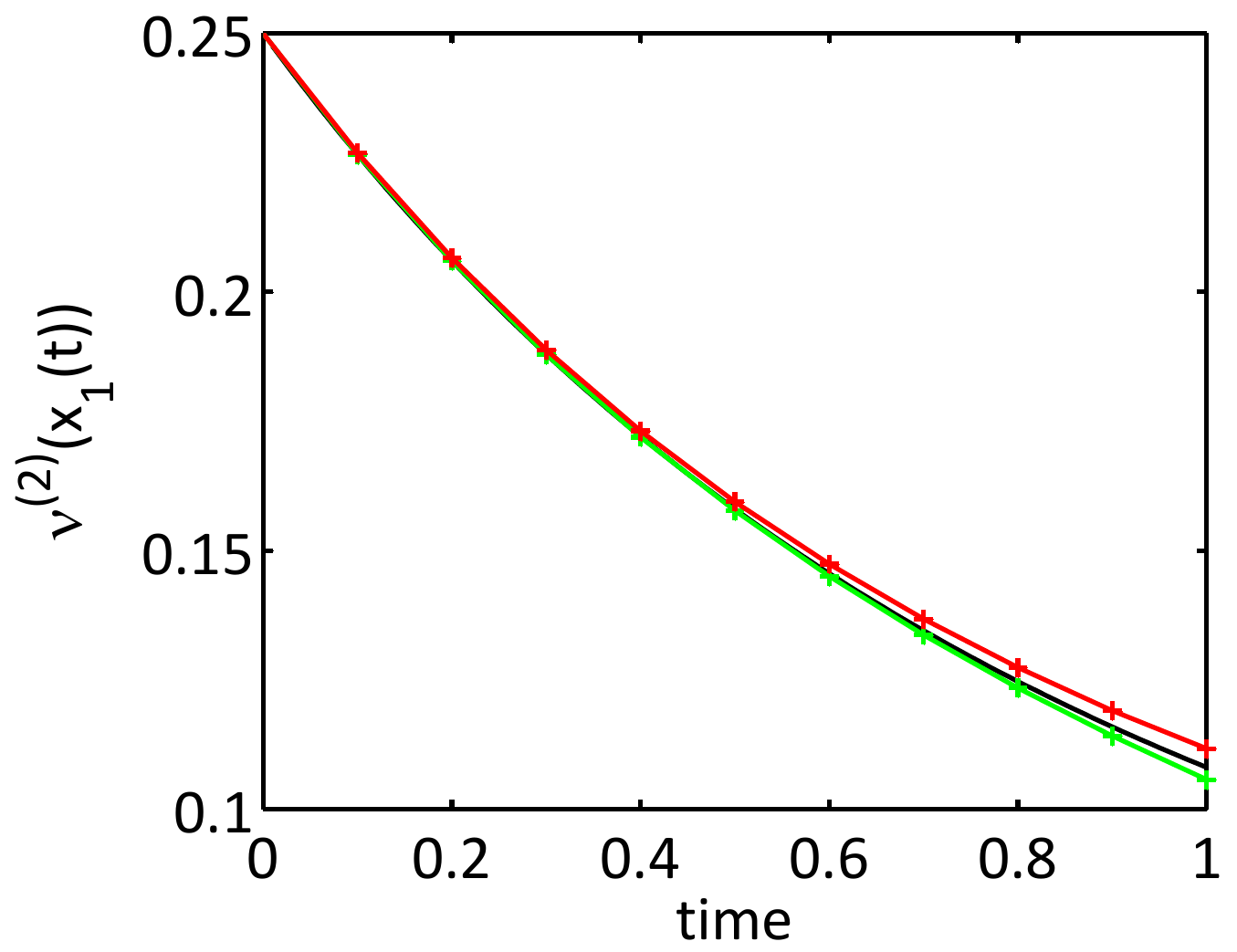}
\caption{Estimation of moments for example 1. The pointwise-in-time estimates of the second moment $\nu^{(2)}$ of $x_2(t)$ and the comparison with the analytic solution (black lines) for relaxation order $d = 3$ are shown. The green and red lines correspond to the lower and upper bound, respectively. The pointwise-in-time estimates were connected by lines to guide the eye. Estimates for the first moment are not shown but the estimates are even tighter.}
\label{fig:mom_estim:example}
\end{figure}

We estimated the moments pointwise-in-time at $t=0,0.1,\ldots, 1.0$. The results are presented in Fig.~\ref{fig:mom_estim:example} and compared it with the analytic solution \eqref{eq:mom_estim:example1:solution}. As can be seen, the moment estimates are very tight already for small relaxation orders.

\subsection{Estimation of Probability Mass}
\label{sec:MInv}

This section addresses Prob.~\ref{prob:MInv}, which is approached similar as the pointwise-in-time estimations of the moments. 
However, instead of splitting the occupation measure in time, the measure of interest (\eg the initial measure $\mu_0$) is split in space, which gives different measures supported on the different subdomains and which are linked by Liouville's equation.
To estimate the lower and upper bounds on the probability masses over the subdomains, the $0^\text{th}$ moment of the corresponding measure is minimized and maximized, respectively.
As above, the infinite-dimensional decision problem can be relaxed by a converging hierarchy of semidefinite programs, \cf Eq.~\eqref{eq:prob1:mom_relax}. The results are illustrated in the following example.

\subsubsection{Example 2.}
\label{sec:MInv:example}

\begin{figure}[h]
\centering
\textbf{(a)}
\includegraphics[width=0.9\linewidth]{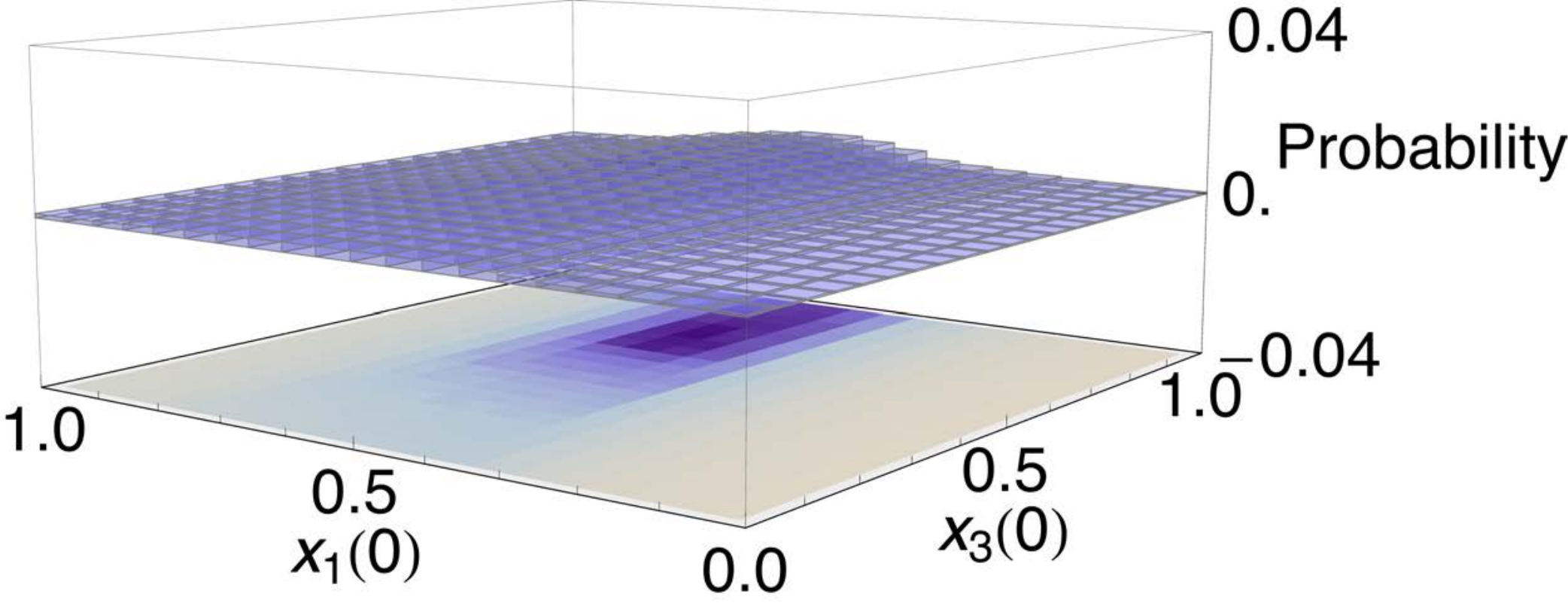}\\
\textbf{(b)}
\includegraphics[width=0.9\linewidth]{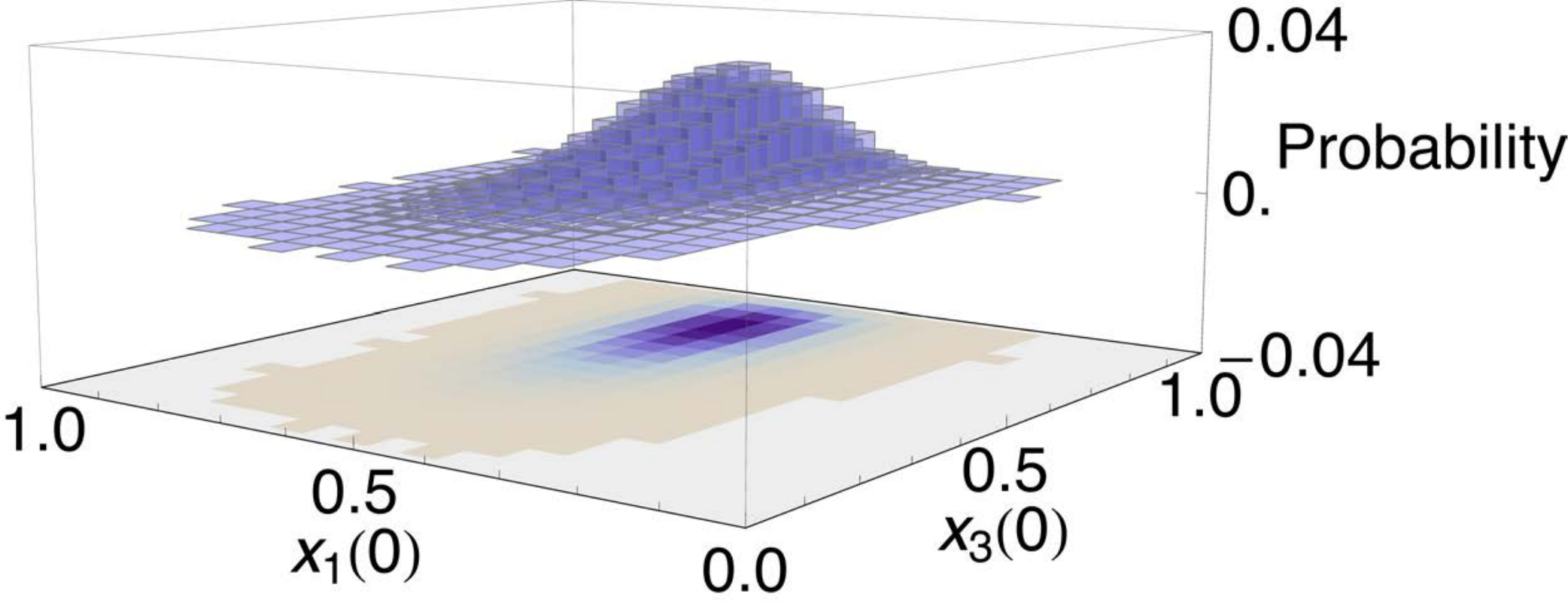}
\caption{Probabilistic model validation. Upper bounds on the probability mass on the indicated box-shaped subsets determined by the proposed approach (a) and the distribution of samples used for the Monte-Carlo simulations for the generation of data (b).}
\label{fig:MInv:example}
\end{figure}

Consider the example 
\begin{subequations}
\begin{align}
\label{eq:mom_estim:example1b}
\dot{x}_1 &= -x_1^3 + x_1x_2 \\
\dot{x}_2 &= -x_2^2 - x_1x_2x_3 + 1 \\
\dot{x}_3 &= 0.
\end{align}
\end{subequations}
The initial condition of $x_2$ was assumed fixed ($x_2(0) = 0.5$), and the initial condition of $x_1$ and of parameter $x_3$ were assumed to be distributed according to Beta distributions (see Fig.~\ref{fig:MInv:example}b):
\begin{align*}
\label{eq:MInv:example:uncertainties}
x_1(0) &\sim \text{Beta}(20,15), \\
x_3(0) &\sim \text{Beta}(5,2).
\end{align*}
Using these distributions we determined the first 10 raw moments for the two states at $t_1 = 0.5$ from 10,000 Monte Carlo simulations and by considering an uncertainty of $1\%$ for the moments. The moments determined from the samples were:
\begin{equation*}
\setlength{\arraycolsep}{1.5pt}
\begin{array}{rclccccccr}
\momv[\le 10]{x_1(0.5)} &=&[ &0.6416 &   0.4156 &   0.2716 &   0.1790 &   0.1188 &\\ &&&   0.0795 &   0.0535 &   0.0362 &   0.0247 &   0.0169 &]\transp, \\
\momv[\le 10]{x_2(0.5)} &=&[ &0.6848 &   0.4695 &   0.3222 &   0.2214  &   0.1523 &\\ &&&   0.1049 &   0.0723 &   0.0499 &   0.0345  &   0.0239 &]\transp.
\end{array}
\end{equation*}
We then partitioned the set $[0,1]\times[0,1]$ into squares of length $1/15$ and determined the lower and upper bound on the probability masses on the subsets. The upper bounds are shown in Fig.~\ref{fig:MInv:example}(a), whereas Fig.~\ref{fig:MInv:example}(b) shows the probabiliy mass of the original distributions. The results demonstrate that the shape of the probability measure as well as the location of its peak can be approximated using the proposed approach.

\section{Conclusion}
\label{sec:conclusion}

This work proposes a combination of set-based and probabilistic uncertainties for the estimation of the support and the moments pointwise-in-time.
In addition, an approach is presented to estimate the probability masses over a subset, which then can be used to approximate the underlying probability distribution and to provide statements on probabilistic model validity.

Due to limitations of state-of-the-art SDP solvers and the high workload of LMIs, only small systems can be treated currently, since the complexity grows polynomially in the number of variables and the relaxation order. However, semi-definite programming and sum-of-squares restrictions is an active research field and computational and algorithmic improvements exists (\eg \cite{Permenter_Parillo_2012_CDC__Selecting_monomial_basis_SOS,Seiler_etAl_2013_Simplification_SOS}) and are likely to be implemented in solvers in the future.

The work provides several directions for future research. For safety and estimation purposes, it might be required to have explicitly time-dependent outer approximations or funnels of the moment and state trajectories. This is particularly challenging for the moments, since these do not explicitly appear in the decision problems.
Another important line of research is the consideration of time-varying disturbances. 
Both extensions are beyond the scope of this article, but interesting future research direction with possible applications for optimal control and estimation of systems governed by stochastic differential equations.

\section*{Acknowledgment}
The authors thank Philipp Rumschinski for very valuable discussions and critical reading of the manuscript, and Milan Korda for helpful discussions. 

                       
\end{document}